\documentclass[12pt]{amsart}
\usepackage{amssymb}
\usepackage{amsfonts}
\newcommand{\nc}{\newcommand}
\nc{\thref}[1]{Theorem~\ref{theo:#1}}
\nc{\selabel}[1]{\label{sect:#1}}
\nc{\seref}[1]{Section~\ref{sect:#1}}
\nc{\lelabel}[1]{\label{lemm:#1}}
\nc{\leref}[1]{Lemma~\ref{lemm:#1}}
\nc{\prlabel}[1]{\label{prop:#1}}
\nc{\prref}[1]{Proposition~\ref{prop:#1}}
\nc{\colabel}[1]{\label{coro:#1}}
\nc{\coref}[1]{Corollary~\ref{coro:#1}}
\nc{\exlabel}[1]{\label{exam:#1}}
\nc{\exref}[1]{Example~\ref{exam:#1}}
\nc{\delabel}[1]{\label{defi:#1}}
\nc{\deref}[1]{Definition~\ref{defi:#1}}
\nc{\eqlabel}[1]{\label{equa:#1}}
\nc{\relabel}[1]{\label{rema:#1}}
\nc{\reref}[1]{Lemma~\ref{rema:#1}}
\providecommand{\operatorname}[1]{\mathrm{#1}\,}
\nc{\Hom}{\operatorname{Hom}} \nc{\Mor}{\operatorname{Mor}}
\nc{\Aut}{\operatorname{Aut}} \nc{\Ann}{\operatorname{Ann}}
\nc{\Ker}{\operatorname{Ker}} \nc{\Trace}{\operatorname{Trace}}
\nc{\Char}{\operatorname{Char}} \nc{\Mod}{\operatorname{Mod}}
\nc{\End}{\operatorname{End}} \nc{\Spec}{\operatorname{Spec}}
\nc{\Span}{\operatorname{Span}} \nc{\sgn}{\operatorname{sgn}}
\nc{\Id}{\operatorname{Id}} \nc{\Com}{\operatorname{Com}}
\nc{\rank}{\operatorname{rank}}

\let\:=\colon

\newtheorem{de}{Definition}[section]
\newtheorem{lm}[de]{Lemma}
\newtheorem{pr}[de]{Proposition}
\newtheorem{co}[de]{Corollary}
\newtheorem{re}[de]{Remark}
\newtheorem{res}[de]{Remarks}
\newtheorem{te}[de]{Theorem}
\newtheorem{ex}[de]{Example}
\newtheorem{exs}[de]{Examples}

\def\bex{\begin{ex}}
\def\eex{\end{ex}}
\def\bexs{\begin{exs}}
\def\eexs{\end{exs}}
\def\bl{\begin{lm}}
\def\el{\end{lm}}
\def\bc{\begin{co}}
\def\ec{\end{co}}
\def\bt{\begin{te}}
\def\et{\end{te}}
\def\bpr{\begin{pr}}
\def\epr{\end{pr}}
\def\br{\begin{re}}
\def\er{\end{re}}
\def\brs{\begin{res}}
\def\ers{\end{res}}
\def\bd{\begin{de}}
\def\ed{\end{de}}
\def\be{\begin{equation}}
\def\ee{\end{equation}}
\def\bea{\begin{eqnarray*}}
\def\eea{\end{eqnarray*}}
\def\bp{\begin{proof}}
\def\ep{\end{proof}}

\let\:=\colon
\begin{document}

\title[\hfil Rational series representations involving $\zeta(2n)$ and Apery's constant]
{Approximations for Apery's constant $\zeta(3)$ and rational series representations involving $\zeta(2n)$}
\begin{abstract}
In this note, using an idea from \cite{Amo-Carrillo-Sanchez} we derive some new series representations involving $\zeta(2n)$ and Euler numbers. Using a well-known series representation for the Clausen function, we also provide some new representations of Apery's constant $\zeta(3)$. In particular cases, we recover some well-known series representations of $\pi$.     

\end{abstract}

\author{Cezar Lupu, Derek Orr}



\maketitle

\section{Introduction and Preliminaries}

In 1734, Leonard Euler produced a sensation when he proved the following formula:

$$\displaystyle\sum_{n=1}^{\infty}\frac{1}{n^2}=\frac{\pi^2}{6}.$$
\vspace{0.3cm}

Later, in 1740, Euler generalized the above formula for even positive integers:

\begin{equation}
\zeta(2n)=(-1)^{n+1}\frac{B_{2n}2^{2n-1}\pi^{2n}}{(2n)!},
\end{equation}

\vspace{0.4cm}

where $\displaystyle\zeta(s)=\sum_{n=1}^{\infty}\frac{1}{n^s}$ is the celebrated Riemann zeta function. The coefficients $B_{n}$ are the so-called Bernoulli numbers and they are defined in the following way:

$$\displaystyle \frac{z}{e^z-1}=\sum_{n=0}^{\infty}\frac{B_{n}}{n!}z^n, |z|<2\pi.$$

\vspace{0.4cm}

An elementary proof of $(1)$ can be found in \cite{Apostol}. In \cite{Amo-Carrillo-Sanchez}, de Amo, Carrillo and Sanchez produced another proof of Euler's formula $(1)$ using the Taylor series expansion of the tangent function and Fubini's theorem. 

In this paper, using similar ideas but for other functions, we provide a new proof for the Clausen acceleration formula that will serve as an application in displaying some fast representations for Apery's constant $\zeta(3)$. Moreover, the Taylor series expansion for the secant and cosecant functions combined with some other integration techniques will give us some interesting rational series representations involving $\zeta(2n)$ and binomial coefficients. Also, we display some particular cases of such series. Some of them are well-known representations of $\pi$. For the sake of completeness we display the Taylor series for the tangent, cotangent, secant and cosecant functions:

\vspace{0.5cm}

\begin{equation}
\tan x=\sum_{n=1}^{\infty}\frac{(-1)^{n-1}2^{2n}(2^{2n}-1)B_{2n}}{(2n)!}x^{2n-1}, |x|<\frac{\pi}{2}
\end{equation}

\begin{equation}
\cot x=\sum_{n=0}^{\infty}\frac{(-1)^n2^{2n}B_{2n}}{(2n)!}x^{2n-1}, |x|<\pi
\end{equation}

\begin{equation}
\sec x=\sum_{n=0}^{\infty}\frac{(-1)^nE_{2n}}{(2n)!}x^{2n}, |x|<\frac{\pi}{2} 
\end{equation}

\begin{equation}
\csc x=\sum_{n=0}^{\infty}\frac{(-1)^{n+1}2(2^{2n-1}-1)B_{2n}}{(2n)!}x^{2n-1}, |x|<\pi
\end{equation}

where $E_{n}$ are the Euler numbers, and $B_{n}$ the Bernoulli numbers.

\vspace{0.4cm}

The Riemann zeta function $\zeta(s)$ and the Hurwitz (generalized) function $\zeta(s, a)$ are defined by

\begin{equation}
\displaystyle\zeta(s):=  \left\{
\begin{array}{ll}
      \displaystyle\sum_{n=1}^{\infty}\frac{1}{n^s}=\frac{1}{1-2^{-s}}\sum_{n=1}^{\infty}\frac{1}{(2n-1)^s} & \operatorname{Re}(s)>1, \\ \\
      \displaystyle\frac{1}{1-2^{1-s}}\sum_{n=1}^{\infty}\frac{(-1)^{n-1}}{n^s}  & \operatorname{Re}(s)>0, s\neq 1,\\
\end{array} 
\right.
\end{equation}

and 

\begin{equation}
\zeta(s, a):=\sum_{n=0}^{\infty}\frac{1}{(n+a)^s}, \operatorname{Re}(s)>1; a\neq 0, -1, -2, \ldots
\end{equation}

\vspace{0.3cm}

Both of them are analytic over the whole complex plane, except $s=1$, where they have a simple pole. Also, from the two definitions above, one can observe that

$$\displaystyle\zeta(s)=\zeta(s, 1)=\frac{1}{2^s-1}\zeta\left(s, \frac{1}{2}\right)=1+\zeta(s, 2).$$

\vspace{0.4cm}

Clausen's function (or Clausen's integral, see \cite{Clausen}) $\displaystyle\operatorname{Cl_{2}}(\theta)$ is defined by

$$\displaystyle\operatorname{Cl_{2}}(\theta):=\sum_{k=1}^{\infty}\frac{\sin k\theta}{k^2}=-\int_0^{\theta}\log\left(2\sin\left(\frac{x}{2}\right)\right) \hspace{3pt} dx.$$

\vspace{0.4cm}

This integral was considered for the first time by Clausen in 1832 (\cite{Clausen}), and it was investigated later by many authors \cite{Bowman, Choi, Choi-Strivastava, Choi-Strivastava1, Choi-Strivastava-Adamchik, Choi-Cho-Strivastava, Grosjean, Koblig, Wood, Wu-Zhang-Liu}. The Clausen functions are very important in mathematical physics. Some well-known properties of the Clausen function include periodicity in the following sense:

$$\displaystyle\operatorname{Cl_{2}}(2k\pi\pm\theta)=\operatorname{Cl_{2}}(\pm\theta)=\pm\operatorname{Cl_{2}}(\theta).$$

Moreover, it is quite clear from the definition that $\displaystyle\operatorname{Cl_{2}}(k\pi)=0$ for $k$ integer. For example, for $k=1$ we deduce

$$\displaystyle\int_0^{\pi}\log\left(2\sin\Big(\frac{x}{2}\Big)\right) \hspace{3pt} dx=0, \int_0^{\frac{\pi}{2}}\log(\sin x) \hspace{3pt} dx=-\frac{\pi}{2}\log 2.$$

\vspace{0.4cm}

By periodicity we have $\displaystyle\operatorname{Cl_{2}}\left(\frac{\pi}{2}\right)=-\operatorname{Cl_{2}}\left(\frac{3\pi}{2}\right)=G$, where $G$ is the Catalan constant defined by

$$\displaystyle G:=\sum_{n=0}^{\infty}\frac{(-1)^n}{(2n+1)^2}\approx 0.9159... $$

More generally, one can express the above integral as the following:

$$\displaystyle\int_0^{\theta}\log(\sin x) \hspace{3pt} dx=-\frac{1}{2}\operatorname{Cl_{2}}(2\theta)-\theta\log 2,$$
\setlength{\belowdisplayskip}{0pt} 
$$\displaystyle\int_0^{\theta}\log(\cos x) \hspace{3pt} dx=-\frac{1}{2}\operatorname{Cl_{2}}(\pi-2\theta)-\theta\log 2,$$
\setlength{\belowdisplayskip}{0pt} 
$$\displaystyle\int_0^{\theta}\log(1+\cos x) \hspace{3pt} dx=2\operatorname{Cl_{2}}(\pi-\theta)-\theta\log 2,$$

and

$$\int_{0}^{\theta} \log(1+\sin x) \hspace{3pt} dx = 2G-2\operatorname{Cl}_{2}\Big(\frac{\pi}{2}+\theta\Big)-\theta\log2.$$

\bigskip

The Dirichlet beta function is defined as

$$\displaystyle\beta(s):=\sum_{n=0}^{\infty}\frac{(-1)^n}{(2n+1)^s}.$$

Alternatively, one can express the beta function in terms of Hurwitz zeta function by the following formula valid in the whole complex $s$-plane,

$$\displaystyle\beta(s)=\frac{1}{4^s}\left(\zeta(s, 1/4)-\zeta(s, 3/4)\right).$$
\bigskip

Clearly, $\beta(2)=G$ (Catalan's constant), $\beta(3)=\frac{\pi^3}{32}$, and $\displaystyle\beta(2n+1)=\frac{(-1)^nE_{2n}\pi^{2n+1}}{4^{n+1}(2n)!}$,
where $E_{n}$ are the Euler numbers mentioned above which are given by the Taylor series $\displaystyle\frac{2}{e^{t}+e^{-t}}=\sum_{n=0}^{\infty}\frac{E_{n}}{n!}t^n$. 

\subsection{Organization of the paper.} This note is organized as follows. In the first subsection (2.1) of the next section, we give a new proof for the Clausen acceleration formula. This formula serves as an ingredient for Theorem 2.2 where we provide some "fast" converging series representations for Apery's constant $\zeta(3)$. On the other hand, in the second subsection, we produce some rational representations of $\zeta(2n)$ and we also provide some interesting series summations as corollaries. 

\bigskip

\textbf{Acknowledgements.} We would like to thank Piotr Hajlasz, Bogdan Ion, George Sparling and William C. Troy for some fruitful conversations which led to an improvement of the present paper.

\section{Main results}
\bigskip

\subsection{Clausen acceleration formula and some representations of $\zeta(3)$}

We provide a new proof for the classical Clausen acceleration formula \cite{Choi-Strivastava1,Lewin}.

\begin{pr} We have the following representation for the Clausen function $\operatorname{Cl_{2}}(\theta)$,

\begin{equation}\displaystyle\frac{\operatorname{Cl_{2}(\theta)}}{\theta}=1-\log |\theta|+\sum_{n=1}^{\infty}\frac{\zeta(2n)}{(2\pi)^{2n}n(2n+1)}\theta^{2n}, |\theta|<2\pi.
\end{equation}

\end{pr}

\bigskip

\textit{Proof.} Integrating by parts the function $xy\cot(xy)$ we have

$$\displaystyle\int_0^{\frac{\pi}{2}}xy\cot(xy)\hspace{3pt}dx =\frac{\pi}{2}\log\left(\sin\left(\frac{\pi y}{2}\right)\right)-\int_0^{\frac{\pi}{2}}\log(\sin(xy))\hspace{3pt}dx$$

$$\displaystyle =\frac{\pi}{2}\log\left(\sin\left(\frac{\pi y}{2}\right)\right)-\frac{1}{2y}\int_0^{\pi y}\log\left(\sin\left(\frac{u}{2}\right)\right)\hspace{3pt}du$$

$$\displaystyle =\frac{\pi}{2}\log\left(\sin\left(\frac{\pi y}{2}\right)\right)+\frac{\pi}{2}\log 2+\frac{1}{2y}\operatorname{Cl_{2}}(\pi y)$$

$$\displaystyle =\frac{\pi}{2}\log\left(\frac{\pi y}{2}\prod_{k=1}^{\infty}\left(1-\frac{y^2}{4k^2}\right)\right)+\frac{\pi}{2}\log 2+\frac{1}{2y}\operatorname{Cl_{2}}(\pi y) $$

$$\displaystyle =\frac{\pi}{2}\log\left(\frac{\pi y}{2}\right)-\frac{\pi}{2}\sum_{k=1}^{\infty}\sum_{n=1}^{\infty}\frac{\left(\frac{y^2}{4k^2}\right)^n}{n}+\frac{\pi}{2}\log 2+\frac{1}{2y}\operatorname{Cl_{2}}(\pi y),$$

where we have used the product formula for sine, $\sin(\pi x)=\pi x\prod_{k=1}^{\infty}\left(1-\frac{x^2}{k^2}\right)$.

On the other hand, using the formula $\theta\cot\theta=1-2\sum_{n=1}^{\infty}\frac{\zeta(2n)}{\pi^{2n}}\theta^{2n}$ we have

$$\displaystyle\int_0^{\frac{\pi}{2}}xy\cot(xy)\hspace{3pt}dx=\frac{\pi}{2}-2\sum_{n=1}^{\infty}\frac{\zeta(2n)\left(\frac{\pi}{2}\right)^{2n+1}}{\pi^{2n}(2n+1)}y^{2n}.$$

Therefore, we obtain

$$\displaystyle \frac{\pi}{2}-2\sum_{n=1}^{\infty}\frac{\zeta(2n)\left(\frac{\pi}{2}\right)^{2n+1}}{\pi^{2n}(2n+1)}y^{2n}=\frac{\pi}{2}\log(\pi y)-\frac{\pi}{2}\sum_{n=1}^{\infty}\frac{\zeta(2n)}{n4^n}y^{2n}+\frac{1}{2y}\operatorname{Cl_{2}}(\pi y).$$
\vspace{0.3cm}

This implies that

$$\displaystyle\operatorname{Cl_{2}}(\pi y)=2y\left(\frac{\pi}{2}-\pi\sum_{n=1}^{\infty}\frac{\zeta(2n)}{(2n+1)4^n}y^{2n}-\frac{\pi}{2}\log(\pi y)+\frac{\pi}{2}\sum_{n=1}^{\infty}\frac{\zeta(2n)}{n4^n}y^{2n}\right),$$
\vspace{0.3cm}

which, after some computations, is equivalent to

$$\displaystyle\operatorname{Cl_{2}}(\pi y)=\pi y-\pi y\log(\pi y)+\pi y\sum_{n=1}^{\infty}\frac{\zeta(2n)(\pi y)^{2n}}{4^nn(2n+1)\pi^{2n}}.$$

Setting $\alpha=\pi y$, we obtain our result.$\square$

\vspace{0.4cm}

\textbf{Remarks.}
In particular case of $ \theta=\frac{\pi}{2}$, using the fact that $\operatorname{Cl}_{2}(\frac{\pi}{2})=G$, we obtain 

\begin{equation}\displaystyle\sum_{n=1}^{\infty}\frac{\zeta(2n)}{n(2n+1)16^n}=\frac{2G}{\pi}-1+\log\left(\frac{\pi}{2}\right).
\end{equation}
\bigskip

In computations, the following accelerated peeled form \cite{Borwein-Bradley-Crandall} is used

\begin{equation}
\frac{\operatorname{Cl_{2}}(\theta)}{\theta}=3-\log\left(|\theta|\left(1-\frac{\theta^2}{4\pi^2}\right)\right)-\frac{2\pi}{\theta}\log\left(\frac{2\pi+\theta}{2\pi-\theta}\right)+\sum_{n=1}^{\infty}\frac{\zeta(2n)-1}{n(2n+1)}\left(\frac{\theta}{2\pi}\right)^{2n}.
\end{equation}

\vspace{0.3cm}

It is well-known that $\zeta(n)-1$ converges to zero rapidly for large values of $n$. In \cite{Wu-Zhang-Liu}, Wu, Zhang and Liu derive the following representation for the Clausen function $\operatorname{Cl}_{2}(\theta)$,

\begin{equation}
\operatorname{Cl}_{2}(\theta)=\theta-\theta\log\left(2\sin\frac{\theta}{2}\right)-\sum_{n=1}^{\infty}\frac{2\zeta(2n)}{(2n+1)(2\pi)^{2n}}\theta^{2n+1}.
\end{equation}
\bigskip

It is interesting to see that integrating the above formula from $0$ to $\pi/{2}$ we have the following representation for $\zeta(3)$ (see \cite{Choi-Strivastava-Adamchik}), 

\begin{equation}
\zeta(3)=\frac{4\pi^2}{35}\left(\frac{1}{2}+\frac{2G}{\pi}-\sum_{n=1}^{\infty}\frac{\zeta(2n)}{(n+1)(2n+1)16^n} \right).
\end{equation}
\bigskip

Also, in \cite{Strivastava-Glasser-Adamchik}, Strivastava, Glasser and Adamchik derive series representations for $\zeta(2n+1)$ by evaluating the integral $\int_0^{\pi/\omega}t^{s-1}\cot tdt$, $s, \omega\geq 2$ integers in two different ways. One of the ways involves the generalized Clausen functions. When they are evaulated in terms of $\zeta(2n+1)$ one obtains the following formula for $\zeta(3)$,

\begin{equation}
\zeta(3)=\frac{2\pi^2}{9}\left(\log 2+2\sum_{n=0}^{\infty}\frac{\zeta(2n)}{(2n+3)4^n}\right).
\end{equation}

\bigskip

Another remarkable result which led to Apery's proof of the irrationality of $\zeta(3)$ is given by the rapidly convergent series

\begin{equation}
\zeta(3)=\frac{5}{2}\sum_{n=1}^{\infty}(-1)^{n-1}\frac{1}{n^3\binom{2n}{n}}.
\end{equation}

\vspace{0.3cm}

Moreover, in \cite{Cvijovic-Klinowski}, Cvijovic and Klinowski derive the following formula

\begin{equation}
\zeta(3)=-\frac{\pi^2}{3}\sum_{n=0}^{\infty}\frac{(2n+5)\zeta(2n)}{(2n+1)(2n+2)(2n+3)2^{2n}}.
\end{equation}

\vspace{0.3cm}

This formula is related to the one found by Ewell \cite{Ewell},

\begin{equation}
\zeta(3)=-\frac{4\pi^2}{7}\sum_{n=0}^{\infty}\frac{\zeta(2n)}{(2n+1)(2n+2)2^{2n}}.
\end{equation}

\vspace{0.3cm}

The following result will provide some new representations for Apery's constant $\zeta(3)$. The main ingredients in the proof of this next result are Clausen acceleration formulae and Fubini's theorem. 
\bigskip

\bt We have the following series representations

\begin{equation}\displaystyle\zeta(3)=\frac{4\pi^2}{35}\left(\frac{3}{2}-\log\left(\frac{\pi}{2}\right)+\sum_{n=1}^{\infty}\frac{\zeta(2n)}{n(n+1)(2n+1)16^n}\right),
\end{equation}

\begin{equation}
\zeta(3)=-\frac{64}{3\pi}\beta(4)+\frac{8\pi^2}{9}\left(\frac{4}{3}-\log\left(\frac{\pi}{2}\right)+3\sum_{n=1}^{\infty}\frac{\zeta(2n)}{n(2n+1)(2n+3)16^n}\right),
\end{equation}

and

\begin{equation}
\zeta(3)=-\frac{64}{3\pi}\beta(4)+\frac{16\pi^2}{27}\left(\frac{1}{2}+\frac{3G}{\pi}-3\sum_{n=1}^{\infty}\frac{\zeta(2n)}{(2n+1)(2n+3)16^n}\right),
\end{equation}
\vspace{0.3cm}

where $G$ is the Catalan constant, and $\beta(s)$ is the Dirichlet beta function. 

\et

\bigskip

\textit{Proof.} For $(17)$, integrating $(8)$ from $0$ to $\pi/2$ and using Fubini's theorem, we have

$$\displaystyle\int_0^{\frac{\pi}{2}}\operatorname{Cl_{2}}(y) \hspace{3pt} dy=\int_0^{\frac{\pi}{2}}(y-y\log y) \hspace{3pt} dy+\int_0^{\frac{\pi}{2}}\sum_{n=1}^{\infty}\frac{\zeta(2n)}{(2\pi)^{2n}n(2n+1)}y^{2n+1} \hspace{3pt} dy $$ 
\setlength{\belowdisplayskip}{0pt} 
$$ = \frac{\pi^2}{8}-\Big(\frac{\pi^2}{8}\log\Big(\frac{\pi}{2}\Big)-\frac{1}{2}\int_{0}^{\frac{\pi}{2}} y \hspace{3pt} dy\Big) + \sum_{n=1}^{\infty} \frac{\zeta(2n)\big(\frac{\pi}{2}\big)^{2n+2}}{(2\pi)^{2n}n(2n+1)(2n+2)}, $$

which is equivalent to

$$\displaystyle\int_0^{\frac{\pi}{2}}\operatorname{Cl_{2}}(y) \hspace{3pt} dy= \frac{\pi^2}{8}\bigg(\frac{3}{2}-\log\Big(\frac{\pi}{2}\Big)+ \sum_{n=1}^{\infty} \frac{\zeta(2n)}{n(n+1)(2n+1)16^n} \bigg).$$

\vspace{0.3cm}

Alternatively, we can integrate the Clausen function using its definition given in the introduction and changing the order of integration as follows:

$$ \int_{0}^{\frac{\pi}{2}}\operatorname{Cl}_{2}(y) \hspace{3pt} dy = -\int_{0}^{\frac{\pi}{2}}\int_{0}^{y} \log\Big(2\sin\Big(\frac{x}{2}\Big)\Big) \hspace{3pt} dx dy = - \int_{0}^{\frac{\pi}{2}}\int_{x}^{\frac{\pi}{2}} \log\Big(2\sin\Big(\frac{x}{2}\Big)\Big) \hspace{3pt} dy dx$$ 
\setlength{\belowdisplayskip}{0pt} 
$$ = -\int_{0}^{\frac{\pi}{2}}\frac{\pi}{2}\log\Big(2\sin\Big(\frac{x}{2}\Big)\Big) \hspace{3pt} dx + \int_{0}^{\frac{\pi}{2}} x\log2 \hspace{3pt} dx + \int_{0}^{\frac{\pi}{2}} x\log\Big(\sin\Big(\frac{x}{2}\Big)\Big) \hspace{3pt} dx.$$

\vspace{0.3cm}

Using the definition of the Clausen function again and, after the substitution $x=2u$, using the identity $ \displaystyle \int_{0}^{\pi/4} u\log(\sin(u)) \hspace{3pt} du = \frac{35}{128}\zeta(3)-\frac{\pi G}{8} - \frac{\pi^2}{32}\log2$ (see \cite{Choi-Strivastava-Adamchik}), we have

$$ \int_{0}^{\frac{\pi}{2}}\operatorname{Cl}_{2}(y) \hspace{3pt} dy = \frac{\pi}{2}\operatorname{Cl}_{2}\Big(\frac{\pi}{2}\Big)+\frac{\pi^2}{8}\log2 + 4\Big(\frac{35}{128}\zeta(3)-\frac{\pi G}{8}-\frac{\pi^2}{32}\log2\Big) = \frac{35}{32}\zeta(3).$$

\vspace{0.5cm}

Setting the two results equal to each other and solving for $\zeta(3)$ gives us our result. For $(18)$, we proceed similarly to the previous method. We will begin with

$$ \int_{0}^{\pi^2/4} \operatorname{Cl}_{2}(\sqrt{y}) \hspace{3pt} dy = \int_{0}^{\pi^2/4} \big(\sqrt{y}-\sqrt{y}\log(\sqrt{y})\big) \hspace{3pt} dy +\int_{0}^{\pi^2/4} \sum_{n=1}^{\infty} \frac{\zeta(2n)y^{n+1/2}}{n(2n+1)(2\pi)^{2n}} \hspace{3pt} dy$$ 
\setlength{\belowdisplayskip}{0pt} 
$$ = \frac{\pi^3}{12}-\Big(\frac{\pi^3}{12}\log\Big(\frac{\pi}{2}\Big)-\frac{1}{3}\int_{0}^{\pi^2/4} \sqrt{y} \hspace{3pt} dy\Big)+\sum_{n=1}^{\infty} \frac{\zeta(2n)\big(\frac{\pi^2}{4}\big)^{n+3/2}}{n(2n+1)(n+3/2)(2\pi)^{2n}}, $$

\vspace{0.3cm}

which is equivalent to 

$$ \int_{0}^{\pi^2/4} \operatorname{Cl}_{2}(\sqrt{y}) \hspace{3pt} dy = \frac{\pi^3}{12}\bigg(\frac{4}{3}-\log\Big(\frac{\pi}{2}\Big)+3\sum_{n=1}^{\infty} \frac{\zeta(2n)}{n(2n+1)(2n+3)16^n}\bigg). $$

\vspace{0.3cm}

On the other hand, using the definition of the Clausen function and changing the order of integration, we see

$$ \int_{0}^{\pi^2/4} \operatorname{Cl}_{2}(\sqrt{y}) \hspace{3pt} dy = -\int_{0}^{\pi^2/4}\int_{0}^{\sqrt{y}} \log\Big(2\sin\Big(\frac{x}{2}\Big)\Big) \hspace{3pt} dx dy $$ 
$$ = -\int_{0}^{\frac{\pi}{2}}\int_{x^2}^{\pi^2/4} \log\Big(2\sin\Big(\frac{x}{2}\Big)\Big) \hspace{3pt} dy dx = \frac{\pi^2}{4}\operatorname{Cl}_{2}\Big(\frac{\pi}{2}\Big)+\int_{0}^{\pi/2} x^2\log\Big(2\sin\Big(\frac{x}{2}\Big)\Big) \hspace{3pt} dx $$ 
$$ = \frac{\pi^2G}{4} + \frac{1}{768}\bigg(72\pi\zeta(3)-192\pi^2G+\psi_3\Big(\frac{1}{4}\Big)-\psi_3\Big(\frac{3}{4}\Big)\bigg), $$

\vspace{0.3cm}

where $\psi_3$ is the trigamma function. Using the identity $\psi_n(z) = (-1)^{n+1}n!\zeta(n+1,z)$ where $\zeta(k,z)$ is the Hurwitz zeta function, and the relationship between the Hurwitz zeta function and the beta function, we have that

$$ \int_{0}^{\pi^2/4} \operatorname{Cl}_{2}(\sqrt{y}) \hspace{3pt} dy = \frac{3\pi}{32}\zeta(3)+2\beta(4).
$$

\vspace{0.3cm}

Setting this result equal to the previous result of the integral and solving for $\zeta(3)$, we see $(18)$ is indeed true. For $(19)$, instead of integrating the Clausen function using $(8)$, we will integrate it using $(11)$. This gives us,

$$ \int_{0}^{\frac{\pi^2}{4}} \operatorname{Cl}_{2}(\sqrt{y}) \hspace{3pt} dy = \int_{0}^{\frac{\pi^2}{4}} \sqrt{y} \hspace{3pt} dy - \int_{0}^{\frac{\pi^2}{4}} \sqrt{y}\log\Big(2\sin\Big(\frac{\sqrt{y}}{2}\Big)\Big) \hspace{3pt} dy $$
$$\displaystyle +2\int_{0}^{\frac{\pi^2}{4}} \sum_{n=1}^{\infty} \frac{\zeta(2n)y^{n+1/2}}{(2n+1)(2\pi)^{2n}} \hspace{3pt} dy$$ 
$$ = \frac{\pi^3}{12}-2\int_{0}^{\frac{\pi}{2}} u^2\log\Big(2\sin\Big(\frac{u}{2}\Big)\Big) \hspace{3pt} du + 2\sum_{n=1}^{\infty} \frac{\zeta(2n)(\frac{\pi^2}{4})^{n+3/2}}{(2n+1)(n+3/2)(2\pi)^{2n}} $$ 
$$ = \frac{\pi^3}{12}-\frac{1}{384}\Big(72\pi\zeta(3)-192\pi^2G+3!4^4\beta(4)\Big)-\frac{\pi^3}{2}\sum_{n=1}^{\infty} \frac{\zeta(2n)}{(2n+1)(2n+3)16^n}. $$

Setting this result equal to the previous Clausen formula yields

$$\Big(\frac{3\pi}{32}+\frac{3\pi}{16}\Big)\zeta(3) = \frac{9\pi}{32}\zeta(3) = \frac{\pi^3}{12} + \frac{\pi^2G}{2}-6\beta(4)-\frac{\pi^3}{2}\sum_{n=1}^{\infty} \frac{\zeta(2n)}{(2n+1)(2n+3)16^n},$$

\vspace{0.3cm}

and from here, $(19)$ follows. $\square$

\textbf{Remark.} To prove $(17)$, one could solve for $G$ in $(9)$ and plug it into $(12)$ and rearrange. Also, to prove $(19)$, one could solve for $\log\big(\frac{\pi}{2}\big)$ in $(9)$ and plug that into $(18)$ and rearrange. Further, if we integrate $(10)$, we arrive at a rapidly converging series representation for $\zeta(3)$, that is

\begin{equation}
\displaystyle\zeta(3)=\frac{2\pi^2}{35}\bigg(9+138\log2-18\log3-50\log5-2\log\pi+2\sum_{n=1}^{\infty}\frac{\zeta(2n)-1}{n(2n+1)(n+1)16^n}\bigg).
\end{equation} 

\bigskip

\subsection{Rational series representations involving $\zeta(2n)$ and binomial coefficients}

As it has been already been highlighted in \cite{Borwein-Bradley-Crandall} one can relate the rational $\zeta$-series with various Dirichlet $L$-series. A rational $\zeta$-series can be accelerated for computational purposes provided that one solves the exact sum

$$\displaystyle\sum_{n=2}^{\infty}\frac{q_{n}}{a^n},$$

where $a=2, 3, 4, \ldots$.

In fact, it has been already highlighted in \cite{Borwein-Bradley-Crandall} we shall call \textit{rational $\zeta$-series} of a real number $x$, the following representation:

$$\displaystyle x=\sum_{n=2}^{\infty}q_{n}\zeta(n, m),$$

where $q_{n}$ is a rational number and $\zeta(n, m)$ is the Hurwitz zeta function. For $m>1$ integer, one has

$$\displaystyle x=\sum_{n=2}^{\infty}q_{n}\left(\zeta(n)-\sum_{j=1}^{m-1}j^{-n}\right).$$
\bigskip 

In the particular case $m=2$, one has the following series representations:

$$\displaystyle 1=\sum_{n=2}^{\infty}(\zeta(n)-1)$$
\setlength{\belowdisplayskip}{0pt} 
$$\displaystyle 1-\gamma=\sum_{n=2}^{\infty}\frac{1}{n}(\zeta(n)-1)$$
\setlength{\belowdisplayskip}{0pt} 
$$\displaystyle\log 2=\sum_{n=2}^{\infty}\frac{1}{n}(\zeta(2n)-1),$$

where $\gamma$ is the Euler-Mascheroni constant. For other rational zeta series representations we recommend \cite{Adamchik-Strivastava, Choi-Strivastava-Adamchik, Strivastava-Glasser-Adamchik}. 
\bigskip

\bt
The following representation is true

\begin{equation}
\displaystyle\sum_{n=1}^{\infty}\frac{\zeta(2n)}{n4^n}\binom{2n}{m}=  \left\{
\begin{array}{ll}
      \displaystyle\frac{1}{m} & m \operatorname{odd}, \\ \\
      \displaystyle\frac{1}{m}\left(2\zeta(m)\left(1-\frac{1}{2^m}\right)-1\right)  & m \operatorname{even}.\\
\end{array} 
\right.
\end{equation}

\et
\bigskip

\textit{Proof}. We start by integrating $xy\operatorname{csc}(xy)$ two different ways. Applying integration by parts, L'Hospital's rule, and properties of logarithm, we find

$$ \displaystyle\int_{0}^{\frac{\pi}{2}} xy\operatorname{csc}(xy) \hspace{3pt} dx = -\frac{\pi}{2} \Bigg(\log\bigg(1+\cos\Big(\frac{\pi y}{2}\Big)\bigg)-\log\bigg(\sin\Big(\frac{\pi y}{2}\Big)\bigg)\Bigg) $$ 
\setlength{\belowdisplayskip}{0pt} 
$$\displaystyle + \int_{0}^{\frac{\pi}{2}} \log(1+\cos(xy)) \hspace{3pt} dx - \int_{0}^{\frac{\pi}{2}} \log(\sin(xy)) \hspace{3pt} dx$$
\setlength{\belowdisplayskip}{0pt} 
$$= -\frac{\pi}{2}\frac{d}{d\alpha}\Big(2\operatorname{Cl}_{2}(\pi-\alpha)+\frac{1}{2}\operatorname{Cl}_{2}(2\alpha)\Big)+\frac{1}{y}\Big(2\operatorname{Cl}_{2}(\pi-\alpha)+\frac{1}{2}\operatorname{Cl}_{2}(2\alpha)\Big),$$

\bigskip

where $ \alpha = \displaystyle \frac{\pi y}{2}$. Applying $(8)$ to the result and simplifying, we get

$$ \displaystyle\int_{0}^{\frac{\pi}{2}} xy\operatorname{csc}(xy) \hspace{3pt} dx = \pi \sum_{n=1}^{\infty}\frac{\zeta(2n)(\frac{\pi}{2})^{2n}(2-y)^{2n}}{n(2\pi)^{2n}} - \frac{\pi}{2} \sum_{n=1}^{\infty}\frac{\zeta(2n)(\pi y)^{2n}}{n(2\pi)^{2n}}-\frac{\pi}{2}$$ 
\setlength{\belowdisplayskip}{0pt} 
$$ \displaystyle +\frac{2\pi}{y}\big(1-\log\Big(\frac{\pi}{2}\Big)-\log(2-y)\big)+ \frac{2}{y}\sum_{n=1}^{\infty}\frac{\zeta(2n)(\frac{\pi}{2})^{2n+1}(2-y)^{2n+1}}{n(2n+1)(2\pi)^{2n}}+\frac{1}{2y}\sum_{n=1}^{\infty}\frac{\zeta(2n)(\pi y)^{2n+1}}{(2\pi)^{2n}}.$$ 

\vspace{0.3cm}

On the other hand, we can apply Fubini's theorem and integrate its power series term by term. Thus, we will have

$$\displaystyle\int_{0}^{\frac{\pi}{2}} xy\operatorname{csc}(xy) \hspace{3pt} dx = \int_{0}^{\frac{\pi}{2}}\sum_{n=0}^{\infty}\frac{(-1)^{n+1}(4^n-2)B_{2n}}{(2n)!}(xy)^{2n} \hspace{3pt} dx$$ 
\setlength{\belowdisplayskip}{0pt} $$=\sum_{n=0}^{\infty}\frac{(-1)^{n+1}(4^n-2)B_{2n}(\frac{\pi}{2})^{2n+1}}{(2n+1)(2n)!}y^{2n}. $$

\vspace{0.5cm}

Setting the two results equal to each other and simplifying more, we see

$$ \displaystyle \sum_{n=1}^{\infty}\frac{\zeta(2n)}{n4^n}\sum_{k=0}^{2n}\binom{2n}{k}\frac{(-1)^k}{2^k} y^k - \sum_{k=1}^{\infty}\frac{\zeta(2k)}{(2k+1)4^k}y^{2k}+\frac{2}{y}(1-\log\pi)+2\sum_{k=1}^{\infty}\frac{1}{k2^k}y^{k-1}$$ 
\setlength{\belowdisplayskip}{0pt} 
$$\displaystyle -\frac{1}{2} + 2\sum_{n=1}^{\infty}\frac{\zeta(2n)}{n(2n+1)4^n}\sum_{k=0}^{2n+1}\frac{(-1)^k}{2^k}\binom{2n+1}{k}y^{k-1}= \sum_{k=0}^{\infty}\frac{(-1)^{k+1}(2^{2k-1}-1)B_{2k}\pi^{2k}}{4^k(2k+1)!}y^{2k}.$$ 

\vspace{0.4cm}

Now we group the coefficients on both sides. The odd powers of $y$ (i.e., $2j-1$ for $j=1,2,3,...$), we find

$$  \sum_{n=1}^{\infty}\frac{\zeta(2n)}{n4^n}\binom{2n}{2j-1}\frac{(-1)^{2j-1}}{2^{2j-1}}+\frac{2}{(2j)2^{2j}}+2\sum_{n=1}^{\infty}\frac{\zeta(2n)}{n(2n+1)4^n}\frac{(-1)^{2j}}{2^{2j}}\binom{2n+1}{2j}=0.$$

\vspace{0.5cm}

Multiplying by $ 2^{2j-1}$, we see

$$\frac{1}{2j} = \sum_{n=1}^{\infty}\frac{\zeta(2n)}{n4^n}\Bigg(\binom{2n}{2j-1}-\frac{1}{2n+1}\binom{2n+1}{2j}\Bigg)= \frac{2j-1}{2j}\sum_{n=1}^{\infty}\frac{\zeta(2n)}{n4^n}\binom{2n}{2j-1}.$$

\vspace{0.5cm}

Setting $m=2j-1$, we arrive at the first part of the theorem. For the even powers of $y$ (i.e., $2j$ for $j=1,2,3,...$), we arrive at the following:

$$ \sum_{n=1}^{\infty}\frac{\zeta(2n)}{n4^n}\frac{(-1)^{2j}}{2^{2j}}\binom{2n}{2j}-\frac{\zeta(2j)}{(2j+1)4^j}+\frac{2}{(2j+1)2^{2j+1}}$$ 
\setlength{\belowdisplayskip}{0pt} 
$$ \displaystyle +2\sum_{n=1}^{\infty}\frac{\zeta(2n)}{n(2n+1)4^n}\frac{(-1)^{2j+1}}{2^{2j+1}}\binom{2n+1}{2j+1} =\frac{(-1)^{j+1}(2^{2j-1}-1)B_{2j}\pi^{2j}}{4^j(2j+1)!}.$$

\vspace{0.5cm}

Multiplying by $4^j$ and using $(1)$ to replace the Bernoulli numbers by $\zeta(2j)$, we have

$$ \sum_{n=1}^{\infty}\frac{\zeta(2n)}{n4^n}\Bigg(\binom{2n}{2j}-\frac{1}{2n+1}\binom{2n+1}{2j+1}\Bigg)=\frac{\zeta(2j)}{2j+1}-\frac{1}{2j+1}+\frac{\zeta(2j)}{(2j+1)}\Big(1-\frac{2}{2^{2j}}\Big).$$

\vspace{0.5cm}

Using the binomial identity $ \displaystyle \binom{2n}{2j}-\frac{1}{2n+1}\binom{2n+1}{2j+1} = \frac{2j}{2j+1}\binom{2n}{2j}$, this simplifies to 

$$ \sum_{n=1}^{\infty}\frac{\zeta(2n)}{n4^n}\binom{2n}{2j} = \frac{1}{2j}\Bigg(2\zeta(2j)\Big(1-\frac{1}{2^{2j}}\Big)\Bigg).$$
\vspace{0.5cm}

Letting $m=2j$ gives the final result of the theorem and thus, the proof is complete. $\square$

\vspace{0.5cm}

\begin{co} (\cite{Tyler-Chernoff, Boros-Moll}) We have 

\begin{equation} 
\sum_{n=1}^{\infty}\frac{\zeta(2n)}{n(2n+1)4^n}=\log\pi-1.
\end{equation}

\end{co}
\bigskip 

\textit{Proof}. This follows immediately by setting the coefficients of $y^{-1}$ equal to each other on both sides. $\square$
\vspace{0.4cm}

\begin{co}
We have the following series representations

\begin{equation}
\sum_{n=1}^{\infty}\frac{\zeta(2n)}{4^n}=\frac{1}{2},
\end{equation}

\begin{equation}
\sum_{n=1}^{\infty}\frac{\zeta(2n)(2n-1)(2n-2)}{4^n}=1,
\end{equation}

\begin{equation}
\sum_{n=1}^{\infty}\frac{\zeta(2n)(2n-1)}{4^n}=\frac{\pi^2}{8}-\frac{1}{2},
\end{equation}

\begin{equation}
\sum_{n=1}^{\infty}\frac{\zeta(2n)n}{4^n}=\frac{\pi^2}{16},
\end{equation}

\vspace{0.5cm}

and

\begin{equation}
\sum_{n=1}^{\infty}\frac{\zeta(2n)n^2}{4^n}=\frac{3\pi^2}{32}.
\end{equation}

\end{co}
\bigskip 

\textit{Proof}. Define $\displaystyle \operatorname{F}(k) := \sum\limits_{n=1}^{\infty} \frac{\zeta(2n)n^k}{4^n}$. Letting $m=1,3,2$, we obtain $(23)$, $(24)$, and $(25)$, respectively. Using $(23)$, note that $(25)$ can be rewritten as $2F(1)-F(0) = \frac{\pi^2}{8} - F(0)$, and so, $(26)$ follows immediately. Finally, using $(23)$ and $(26)$, note that $(24)$ can be rewritten as $4F(2)-6F(1)+2F(0) = 2F(0)$. From this, $(27)$ follows immediately. $\square$

\textbf{Remark.} Letting $m=2k$ and $m=2k-1$ $(k=1,2,3,...)$ in the even and odd parts of $(21)$ respectively, we can add the two results to find

$$ \sum_{n=1}^{\infty} \frac{\zeta(2n)}{n4^n}\Bigg(\binom{2n}{2k-1}+\binom{2n}{2k}\Bigg) = \frac{\zeta(2k)}{k}\Big(1-\frac{1}{4^k}\Big)+\frac{1}{2k-1}-\frac{1}{2k}, $$

\vspace{0.5cm}

which is equivalent to

\begin{equation}\sum_{n=1}^{\infty}\frac{\zeta(2n)}{n4^n}\binom{2n+1}{2k}=\frac{\zeta(2k)}{k}\left(1-\frac{1}{4^k}\right)-\frac{1}{2k(2k-1)}.
\end{equation}

\bigskip

\bt
We have the following series representation

\begin{equation}
\sum_{n=1}^{\infty}\frac{\zeta(2n)}{n16^n}\binom{2n}{m}=\left\{
\begin{array}{ll}
      \displaystyle \frac{1}{m}\left(1-\zeta_{E}(m-1)\left(1-\frac{1}{2^{m-1}}\right)\right) & m \operatorname{odd}, \\ \\
      \displaystyle \frac{1}{m}\left(\zeta(m)\left(1-\frac{1}{2^m}\right)-1\right)  & m \operatorname{even},\\
\end{array} 
\right.
\end{equation}
\bigskip

where $\displaystyle \zeta_{E}(2k)=\frac{(-1)^{k+1}E_{2k}\pi^{2k+1}}{4(1-4^k)(2k)!}$ and $E_{2k}$ are the Euler numbers.

\et
\bigskip

\textit{Proof}. Similar to the previous theorem's proof, we integrate $\operatorname{sec}(xy)$ two different ways. First, integrating regularly and using properties of logarithm, we have

$$ \int_{0}^{\frac{\pi}{2}} \operatorname{sec}(xy) \hspace{3pt} dx = \frac{1}{y}\log\bigg(1+\sin\Big(\frac{\pi y}{2}\Big)\bigg) - \frac{1}{y}\log\bigg(\cos\Big(\frac{\pi y}{2}\Big)\bigg)$$
$$= \frac{1}{y}\frac{d}{d\alpha}\bigg(-2\operatorname{Cl}_{2}\Big(\frac{\pi}{2}+\alpha\Big)-\alpha\log2\bigg)-\frac{1}{y}\log\Bigg(\prod_{n=1}^{\infty}\bigg(1-\Big(\frac{y}{2n-1}\Big)^2\bigg)\Bigg),$$

\bigskip

where $\displaystyle \alpha = \frac{\pi y}{2}$ and we have used the product formula $\cos(\beta) =  \prod\limits_{k=1}^{\infty}\Big(1-\big(\frac{2\beta}{\pi(2k-1)}\big)^2\Big)$. Using equation $(8)$ in the first term and applying properties of logarithm and its power series to the second term, the integral becomes

$$ \displaystyle \int_{0}^{\frac{\pi}{2}} \operatorname{sec}(xy) \hspace{3pt} dx = \frac{1}{y}\bigg(2\log\Big(\frac{\pi}{2}\Big)+2\log(1+y)-\log2\bigg) $$
\setlength{\belowdisplayskip}{0pt} 
$$ \displaystyle -\frac{2}{y}\sum_{n=1}^{\infty} \frac{\zeta(2n)(\frac{\pi}{2})^{2n}(1+y)^{2n}}{n(2\pi)^{2n}} \displaystyle +\sum_{n=1}^{\infty}\sum_{k=1}^{\infty}\frac{y^{2k}}{k(2n-1)^{2k}}, $$

\vspace{0.5cm}

and 

$$\displaystyle \sum\limits_{n=1}^{\infty} \frac{1}{(2n-1)^{2k}} = \sum\limits_{n=1}^{\infty} \frac{1}{n^{2k}} - \sum\limits_{n=1}^{\infty} \frac{1}{(2n)^{2k}} = \Big(1-\frac{1}{4^k}\Big)\zeta(2k).$$ 

\bigskip

Alternatively, we can apply Fubini's theorem once again and integrate the power series of the secant function term by term. Doing so will give us the following:

$$\displaystyle\int_{0}^{\frac{\pi}{2}} xy\operatorname{sec}(xy) \hspace{3pt} dx = \int_{0}^{\frac{\pi}{2}}\sum_{n=0}^{\infty}\frac{(-1)^nE_{2n}}{(2n)!}(xy)^{2n} \hspace{3pt} dx = \sum_{n=0}^{\infty}\frac{(-1)^nE_{2n}(\frac{\pi}{2})^{2n+1}}{(2n+1)(2n)!}y^{2n}. $$

\bigskip

Setting the two results equal to each other and simplifying more, we see

$$ \displaystyle \frac{2}{y}\log\Big(\frac{\pi}{2\sqrt{2}}\Big)+2\sum_{k=1}^{\infty}\frac{(-1)^{k+1}}{k}y^{k-1}-2\sum_{n=1}^{\infty}\frac{\zeta(2n)}{n16^n}\sum_{k=0}^{2n}\binom{2n}{k}y^{k-1} $$ 
\setlength{\belowdisplayskip}{0pt} 
$$ \displaystyle +\sum_{k=1}^{\infty}\frac{\zeta(2k)}{k}\bigg(1-\frac{1}{4^k}\bigg)y^{2k-1} = \sum_{k=0}^{\infty}\frac{(-1)^kE_{2k}(\frac{\pi}{2})^{2k+1}}{(2k+1)!}y^{2k}.$$

\bigskip

Now we can group coefficients. For the even powers of y (i.e., $2j$ for $j=1,2,3,...$), we have

$$ 2\frac{(-1)^{2j+2}}{2j+1} - 2\sum_{n=1}^{\infty} \frac{\zeta(2n)}{n16^n}\binom{2n}{2j+1} = \frac{(-1)^jE_{2j}\pi^{2j+1}}{(2j+1)!2^{2j+1}}.$$

\bigskip

Define $\displaystyle \zeta_{E}(2j) := \frac{(-1)^{j+1}E_{2j}\pi^{2j+1}}{4(1-4^j)(2j)!}.$ The motivation comes from Euler's formula for $\zeta(2n)$ and the asymptotic formula $ \displaystyle B_{2n} \sim \frac{\pi E_{2n}}{2^{2n+1}(1-4^n)}, n\to\infty$. For future computations, note that $ \displaystyle \zeta_{E}(2) = \frac{\pi^3}{24}$. Now, simplifying the above expression,

$$\sum_{n=1}^{\infty} \frac{\zeta(2n)}{n16^n}\binom{2n}{2j+1} = \frac{1}{2j+1} + \frac{4\zeta_{E}(2j)(1-4^j)(2j)!}{4(2j+1)!4^j} = \frac{1}{2j+1}\bigg(1-\zeta_{E}(2j)\Big(1-\frac{1}{4^j}\Big)\bigg). $$

\bigskip

Setting $m=2j+1$, we achieve the first result of the theorem. For the odd powers of y (i.e., $2j-1$ for $j=1,2,3,...$), we find

$$ 2\frac{(-1)^{2j}}{2j}-2\sum_{n=1}^{\infty} \frac{\zeta(2n)}{n16^n}\binom{2n}{2j}+\frac{\zeta(2j)}{j}\bigg(1-\frac{1}{4^j}\bigg) = 0,$$
\vspace{0.3cm}

and rearranging this gives the second part of the theorem. $\square$

\bigskip

\begin{co}
We have the following representations

\begin{equation}
\sum_{n=1}^{\infty}\frac{\zeta(2n)}{n16^n}=\log\left(\frac{\pi}{2\sqrt{2}}\right),
\end{equation}

\vspace{0.3cm}

and

\begin{equation}
\sum_{n=1}^{\infty}\frac{\zeta(2n)}{16^n}=\frac{4-\pi}{8}.
\end{equation}

\end{co}

\vspace{0.3cm}

\textit{Proof}.  Formula $(30)$ follows immediately by setting the coefficients of $y^{-1}$ equal to each other on both sides. For $(31)$, we set the constant terms on both sides equal to find 

$$ 2-4\sum\limits_{n=1}^{\infty} \frac{\zeta(2n)}{16^n} = \frac{\pi}{2}.$$

\vspace{0.5cm}

From here, $(31)$ follows immediately. One could also use the theorem for $m=2k+1$ for $k=0$ and note that $\lim_{k\to0} \zeta_{E}(2k)(1-4^k) = \pi/4$. $\square$

\textbf{Remark.} In the previous theorem, if we used the Clausen identity for $\log(\cos(x))$ given in the introduction rather than the cosine product formula, we obtain the following relation:

$$ \log\Big(\frac{\pi}{4}\Big)-\sum_{n=1}^{\infty} \frac{\zeta(2n)}{n4^n}\sum_{k=0}^{2n}\binom{2n}{k}y^k\Big(\frac{2}{4^n}-(-1)^k\Big) +\sum_{k=1}^{\infty} \frac{y^k}{k}(1-2(-1)^k)$$
\setlength{\belowdisplayskip}{0pt}
$$= \sum_{k=0}^{\infty}\frac{2\zeta_{E}(2k)}{2k+1}\Big(1-\frac{1}{4^k}\Big)y^{2k+1}.$$

\bigskip

Gathering the coefficients of $y^0$, we find

$$ \log\Big(\frac{\pi}{4}\Big)-\sum_{n=1}^{\infty}\frac{\zeta(2n)}{n4^n}\Big(\frac{2}{4^n}-1\Big)=0.$$

\bigskip

Rearranging and using $(30)$, we can see that

\begin{equation}
\sum_{n=1}^{\infty}\frac{\zeta(2n)}{n4^n}=\log\left(\frac{\pi}{2}\right),
\end{equation}

\vspace{0.5cm}

which is a very nice result. This series appears by setting the coefficients of $y^0$ equal to each other in the proof of Theorem 2.3; however, the series will cancel and you arrive at a truth statement. Here, we are able to recover that series.

\begin{co}
We have

\begin{equation}
\sum_{n=1}^{\infty}\frac{\zeta(2n)(2n-1)}{16^n}=\frac{\pi^2}{16}-\frac{1}{2},
\end{equation}

\begin{equation}
\sum_{n=1}^{\infty}\frac{\zeta(2n)(2n-1)(2n-2)}{16^n}=1-\frac{\pi^3}{96},
\end{equation}

\begin{equation}
\sum_{n=1}^{\infty}\frac{\zeta(2n)n}{16^n}=\frac{\pi}{16}\left(\frac{\pi}{2}-1\right),
\end{equation}
\vspace{0.5cm}

and

\begin{equation}
\sum_{n=1}^{\infty}\frac{\zeta(2n)n^2}{16^n}=\frac{\pi}{32}\left(\frac{3\pi}{2}-\frac{\pi^2}{4}-1\right).
\end{equation}

\end{co}

\textit{Proof}. Define $\displaystyle \operatorname{G}(k) := \sum\limits_{n=1}^{\infty} \frac{\zeta(2n)n^k}{16^n}$. Letting $m=2,3$ in $(29)$ gives us formulas $(33)$ and $(34)$, respectively. Using $(31)$, we can rewrite $(33)$ as $\displaystyle 2G(1)-G(0)=\frac{\pi^2}{16}-G(0)-\frac{\pi}{8}$ and from this, $(35)$ follows immediately. Lastly, for $(36)$, expanding $(34)$ using $(31)$, we have $\displaystyle 4G(2)-6G(1)+2G(0) = \frac{\pi}{4}+2G(0)-\frac{\pi^3}{96}$. Using $(35)$, we achieve the desired result. $\square$

\bigskip

\begin{co}
We have the following series

\begin{equation}
\sum_{n=1}^{\infty}\frac{\zeta(2n)}{n4^n}\left(1-\frac{1}{4^n}\right)\binom{2n}{2k}=\frac{\zeta(2k)}{2k}\left(1-\frac{1}{4^k}\right),
\end{equation}

\vspace{0.5cm}

and

\begin{equation}
\sum_{n=1}^{\infty}\frac{\zeta(2n)}{n4^n}\left(1-\frac{1}{4^n}\right)\binom{2n}{2k+1}=\frac{\zeta_{E}(2k)}{2k+1}\left(1-\frac{1}{4^k}\right).
\end{equation}

\end{co}
\vspace{0.3cm}

\textit{Proof}. Denote $(29.1)$ and $(29.2)$ as well as $(21.1)$ and $(21.2)$ as the formula for $m$ odd and $m$ even, respectively. Letting $m=2k$ for $k=1,2,3...$ and subtracting $(29.2)$ from $(21.2)$, formula $(37)$ follows immediately. Similarly, letting $m=2k+1$ for $k=0,1,2,...$ and subtracting $(29.1)$ from $(21.1)$, the obtain $(38)$. $\square$ 

\textbf{Remark.} You can add $(37)$ and $(38)$ together to get $\displaystyle \binom{2n+1}{2k+1}$ inside the series. You can also change $k$ to $k-1$ in $(38)$ and add $(37)$ and $(38)$ again to give you $\displaystyle \binom{2n+1}{2k}$ in the series.

\bigskip

\vspace*{3mm}
\begin{flushright}
\begin{minipage}{148mm}\sc\footnotesize
University of Pittsburgh, Department of Mathematics, 301 Thackeray Hall, Pittsburgh, PA 15260, USA\\
{\it E--mail address}: {\tt cel47@pitt.edu, lupucezar@gmail.com} \vspace*{3mm}
\end{minipage}
\end{flushright}

\begin{flushright}
\begin{minipage}{148mm}\sc\footnotesize
University of Pittsburgh, Department of Mathematics, 301 Thackeray Hall, Pittsburgh, PA 15260, USA\\
{\it E--mail address}: {\tt djo15@pitt.edu} \vspace*{3mm}
\end{minipage}
\end{flushright}


\begin{thebibliography}{199}
\markboth{Bibliography}{Bibliography}
%\markboth{References}{References}

%%%%%%%%%%%%%%%%%%%%%%%%%%%%%%%%%%%%%%%%%%%%%%%%%%%%%%%%%%%%%%%

\bibitem{Abramowitz-Stegun}
M. Abramowitz, I. A. Stegun, {\em Handbook of Mathematical Functions with Formulas, Graphs and Mathematical Tables}, New York, NY:Dover, 1972.

\bibitem{Adamchik}
V. S. Adamchik, Contributions to the theory of the Barnes function, {\em International J. Math. Comp. Science} \textbf{9} (2014), 11-30.

\bibitem{Adamchik-Strivastava}
V. S. Adamchik, H. M. Strivastava, Some series of the zeta and related functions, {\em Analysis}, \textbf{18} (1998), 131--144.

\bibitem{Amo-Carrillo-Sanchez}
E. De Amo, M. Diaz Carrillo, J. Hernandez-Sanchez, Another proof of Euler's formula for $\zeta(2k)$, {\em Proc. Amer. Math. Soc.} \textbf{139} (2011), 1441--1444.

\bibitem{Apostol}
T. M. Apostol, Another elementary proof of Euler's formula for $\zeta(2n)$, {\em Amer. math. Monthly}, textbf{80} (1973), 425--431.


\bibitem{Boros-Moll}
G. Boros, V. Moll, {\em Irresistible Integrals. Symbolics, analysis and experiments in the evaluation of integrals}, Cambridge University Press, Cambridge, 2004.

\bibitem{Borwein-Bradley-Crandall}
J. M. Borwein, D. M. Bradley, R. E. Crandall, Computational strategies for the Riemann zeta function, {\em J. Comp. Appl. Math.} \textbf{121} (2000), 247--296.

\bibitem{Borwein-Broadhurst-Kamnitzer}
J. M. Borwein, D. J. Broadhurst, J. Kamnitzer, Central binomial sums, multiple Clausen values, and zeta values, {\em Exp. Math.} \textbf{10} (2001), 25--34. 

\bibitem{Bowman}
F. Bowman, Note on the integral $\int_0^{\frac{\pi}{2}}(\log\sin\theta)^nd\theta$, {\em J. London Math. Soc.} \textbf{22} (1947), 172--173.

\bibitem{Choi}
J. Choi, Some integral representations of the Clausen function $\operatorname{Cl_{2}}(x)$ and the Catalan constant, {\em East Asian Math. J.} \textbf{32} (2016), 43--46.

\bibitem{Choi-Strivastava}
J. Choi, H. Strivastava, Certain classes of series involving the zeta function, {\em J. Math. Anal. Appl.} \textbf{231} (1999), 91--117.

\bibitem{Choi-Strivastava1}
J. Choi, H. M. Strivastava, The Clausen function $\operatorname{Cl_{2}}(x)$ and its related integrals, {\em Thai J. Math.} \textbf{12} (2014), 251--264.

\bibitem{Choi-Strivastava-Adamchik}
J. Choi, H. M. Strivastava, V. S. Adamchik, Multiple Gamma and related functions, {\em Appl. Math. Comput.} \textbf{134} (2003), 515--533. 

\bibitem{Choi-Cho-Strivastava}
J. Choi, Y. J. Cho, H. M. Strivastava, Log-sine integrals involving series associated with the zeta function and polylogarithm function, {\em Math. Scand.} \textbf{105} (2009), 199--217. 

\bibitem{Clausen}
T. Clausen, Uber die function $\sin\phi+\frac{1}{2^2}\sin 2\phi+\frac{1}{3^2}\sin 3\phi+etc.$, {\em J. Reine Angew. Math.} \textbf{8} (1832), 298--300.

\bibitem{Coffey}
M. W. Coffey, One one-dimensional digamma and polygamma series related to the evaluation of Feynman diagrams, {\em J. Comp. Appl. Math.} \textbf{183} (2005), 84--100.

\bibitem{Cvijovic}
D. Cvijovic, New integral representations of the polylogarithm function, {\em Proc. Royal Soc. A} \textbf{643} (2007), 897--905.

\bibitem{Cvijovic1}
D. Cvijovic, Closed-form evaluation of some families of cotangent and cosecant integrals, {\em Integral Transforms Spec. Func.} \textbf{19} (2008), 147--155.

\bibitem{Cvijovic-Klinowski}
D. Cvijovic, J. Klinowski, New rapidly convergent series representations for $\zeta(2n+1)$, {\em Proc. Amer. Math. Soc.} \textbf{125} (1997), 1263--1271.

\bibitem{Doeder}
P. J. de Doeder, On the Clausen integral $\operatorname{Cl_{2}}(\theta)$ and a related integral, {\em J. Comp. Appl. Math.} \textbf{11} (1984), 325--330.

\bibitem{Ewell}
J. A. Ewell, A new series representation for $\zeta(3)$, {\em Amer. Math. Monthly} \textbf{97} (1990), 219--220.

\bibitem{Grosjean}
C. C. Grosjean, Formulae concerning the computation of the Clausen integral $\operatorname{Cl_{2}}(\theta)$, {\em J. Comp. Appl. Math.} \textbf{11} (1984), 331--342.

\bibitem{Koblig}
K. S. Koblig, Chebyshev coefficients for the Clausen function $\operatorname{Cl_{2}}(x)$, {\em J. Comput. Appl. Math.} \textbf{64} (1995), 295--297.

\bibitem{Lewin}
L. Lewin, {\em Polylogarithms and Associated functions}, Elsevier North Holland Inc. New York, 1981.

\bibitem{Lupu}
C. Lupu, Multiple zeta values and infinite series of Leibniz type involving $\zeta(2n)$, {\em in preparation}.

\bibitem{Ogreid-Osland}
O. M. Ogreid, P. Osland, Some infinite series related to Feynman diagrams, \textbf{140} (2002), 659--671.

\bibitem{Strivastava-Glasser-Adamchik}
H. M. Strivastava, M. L. Glasser, V. S. Adamchik, Some definite integrals associated with the Riemann zeta function, {\em Z. Anal. Anwendungen} \textbf{19} (2000), 831--846.

\bibitem{Tyler-Chernoff}
D. Tyler, P. R. Chernoff, An old sum reappears-Elementary problem 3103, {\em Amer. Math. Monthly} \textbf{92} (1985), 507.

\bibitem{Wood}
V. E. Wood, Efficient calculation of Clausen's integral, {\em Math. Comp.} \textbf{22} (1968), 883--884.

\bibitem{Wu-Zhang-Liu}
J. Wu, X. Zhang, D. Liu, An efficient calculation of the Clausen functions, {\em BIT Numer. Math. } \textbf{50} (2010), 193--206.



\end{thebibliography}
\end{document}